\documentclass[12pt]{amsart}
\renewcommand{\bar}{\overline}

\newcommand{\pa}{\partial}

\usepackage{amsmath,amsthm,amscd}

\newcommand{\frk}[1]{{\mathfrak{#1}}}

\title[]{On the Geometry of Classifying Spaces
and Horizontal Slices} 
\author[]{Zhiqin Lu} 
\date{Nov. 14, 1997}
\address[Zhiqin Lu]
{Department of Mathematics\\
Columbia University\\
New York, NY 10027}
\email{lu@math.columbia.edu}

\newtheorem{theorem}{Theorem}[section]
\newtheorem{lemma}{Lemma}[section]
\newtheorem{cor}{Corollary}[section]
\newtheorem{prop}{Proposition}[section]

\newtheorem{definition}{Definition}[section]

\theoremstyle{remark}
\newtheorem{rem}{Remark}[section]
\begin{document}
\maketitle
\tableofcontents
\renewcommand{\baselinestretch}{1.5}

\numberwithin{equation}{section}

\section{Introductions}
Let $(X,\omega)$ be a polarized simply connected  Calabi-Yau
manifold. That is, $X$ is a simply connected compact K\"ahler manifold of
dimension $n$ with zero first Chern class and $\omega$ is a K\"ahler
form
of $X$ such that $[\omega]\in H^2(X,Z)$. In this paper, we study the 
local properties of the moduli
space ${\mathcal M}$ of the polarized Calabi-Yau manifold $(X,\omega)$. 
By definition ${\mathcal M}$ is the parameter space of the complex
structures over $X$ for the fixed polarization $[\omega]$. 
${\mathcal M}$ is a quasi-projective variety by a theorem of
Viehweg~\cite{V2}

Suppose  $X\in
{\mathcal M}$ is the Calabi-Yau manifold. Let
\[
N=\dim\{\eta\in H^1(X,T_X)|\eta\lrcorner \omega=0\}
\]
where $T_X$ is the holomorphic tangent bundle of $X$. By a theorem of
Tian~\cite{T1}, we know that the polarized  universal deformation space
of the
Calabi-Yau
manifold is smooth and has dimension $N$. Since for each $X'\in{\mathcal
M}$, there are no nonzero holomorphic tangent vectors on $X'$, we
concluded
that in each neighborhood 
$\tilde U$of $X'$, there is an open neighborhood $U$ in
$C^N$ such that $U$ is a finite covering of $\tilde U$. Thus the moduli
space is a complex orbifold. A good reference of the  
theorem of Tian can be found in~\cite{Fd}.

By the above observation, we need only to study the local properties of
$U$,
since $U$ is only a finite covering of $\tilde U$, it will be easy to pass
the properties of $U$ to $\tilde U$.
Let $D$ be the
classifying space
corresponding to the Calabi-Yau manifold $X$. The period map
$U\rightarrow D$ is a holomorphic immersion. Furthermore, 
the image of $U$  under the period map is a
horizontal slice of $D$. That is, $U$ is a complex integral submanifold
of the horizontal distribution of $D$. For the precise definitions
of the classifying space and the horizontal slices, see
section 2-3.

Theorem~\ref{thm1} and \ref{thm2} are given below for the horizontal
slices. Thus our results are more general than merely on  moduli
spaces. However, few examples of horizontal slices can be found which are
not moduli spaces because of the non-integrability of the horizontal
distributions.

Let $D=G/V$ where $G$ is a semi-simple Lie group of noncompact type
without compact factors,  and
$V$ is a  connected compact subgroup. Let $K$ be the maximum connected
compact subgroup containing $V$. Then $D_1=G/K$ is a symmetric space.

In this paper, we proved the following:

\begin{theorem}\label{thm1}
Suppose $U\subset D$ is a horizontal slice. Then
\[
p: U\rightarrow D_1
\]
defined by the  natural projection $D\rightarrow
D_1$ is again an
immersion. Furthermore, it is a pluriharmonic map. i.e. It restricts to a
harmonic map from any holomorphic curve of $U$. 
Or in other word, it satisfies
\[
\nabla_{p_*X}{p_*X}
+\nabla_{p_*JX}{p_*JX}
+p_*J[X,JX]=0
\]
for any  vector field $X\in T(U)$, where $J$ is the complex 
structure on $U$ and $\nabla$ is the 
Riemannian connection on $G/K=D_1$.
\end{theorem}

Theorem ~\ref{thm1} is of its own interest. But the most important
application of theorem~\ref{thm1} is the following:

\begin{theorem}\label{thm3}
The restriction of the invariant Riemannian metric on $D_1$ to $U$ is
K\"ahler.
Moreover, the holomorphic bisectional curvature of such a metric
is nonpositive. Furthermore, the Ricci curvature 
and the holomorphic sectional curvature are negative away from
zero.
\end{theorem}

In the last section of this paper, we proved
the non-existence of invariant K\"ahler
metric on some classifying spaces.

Suppose $\Gamma$ is a lattice of the group $G$. We call that $\Gamma$ is
cocompact, if $\Gamma\backslash G$ is a compact topological space.

\begin{theorem}\label{thm2}
Let $D$ be the classifying space of some Calabi-Yau threefold. Then
there are no K\"ahler metrics on $D$ which are $\Gamma$ invariant if
$\Gamma$ is cocompact.
\end{theorem}

The materials are arranged as follows: In section 2-3 we
set up the definitions and notations and
 proved some basic
properties of  horizontal slices,  classifying space and the period
map. Theorem~\ref{thm1} was proved in Section 4. In section 5, we proved
that the Hodge metric
is K\"ahler and the Ricci curvature is negative away from zero,
which is Theorem~\ref{thm3}. In the
last section, we proved Theorem~\ref{thm2} of  the non-existence of the
K\"ahler metric with
some kind of invariance.

The idea to prove theorem ~\ref{thm1} and theorem~\ref{thm2} is to
consider
the
K\"ahler form $\omega$ on the classifying space $D$. 
 It is known that {\it if} the projection $D\rightarrow D_1$
is
holomorphic, then $D$ is K\"ahlerian and $\omega$ is the K\"ahler
form. However, $D\rightarrow D_1$ is in general not
holomorphic. Nevertheless, there is a relation between $\omega$ and
the pull back of the invariant metric on $D_1$. From this we
concluded that although $d\omega\neq 0$ as a
differential form, at some directions, $d\omega$ is indeed
zero. In particular, we use the fact that if $X,Y,Z$ are horizontal, then
$d\omega(X,Y,Z)=0$ in Theorem~\ref{thm1}, and we use the fact that if
$X,Y$
are vertical and $Z$ is horizontal, then $d\omega(X,Y,Z)=0$ 
in Theorem~\ref{thm2}.

The similar result to  Theorem~\ref{thm1} of the 
pluriharmonicity 
of the projection was studied in Bryant~\cite{Br}, Burstall
and Salamon~\cite{BS}, Black~\cite{Bl}. Those papers only considered the
compact cases and can not apply to the cases we are interested in this
paper.

{\bf Acknowledgment} This paper is a refinement of a part of my ph.\,D
thesis. The author would like to thank his advisor, Professor
G. Tian for his advise and constant encouragement during my four
year's ph.\,D study. He also thanks Professor S. T. Yau for his constant
encouragement and many important ideas stimulating further studies of this
problem. He also thanks Professor J. Jost for 
the paper~\cite{JY1}  and pointing out some mistakes
in the references.

\section{The Hodge Structure and the Classifying Space}

Let $X$ be a compact K\"ahler manifold
of dimension $n$. A $C^\infty$ form on $X$ 
decomposes into $(p,q)$-components according to the number of 
$dz's$ and $d\bar{z}'s$. Denoting the $C^\infty$ $n$-forms and
the 
$C^\infty(p,q)$ forms on $X$ by $A^n(X)$ and $A^{p,q}(X)$ 
respectively, we have the decomposition
\[
A^n(X)=\underset{p+q=n}{\oplus}A^{p,q}(X)
\]

The cohomology group is defined as
\begin{align*}
H^{p,q}(X)=&\{closed (p,q) -forms\}/\{exact (p,q)-forms\}\\
=&\{\phi\in A^{p,q}(X)|d\phi=0\}/dA^{n-1}(X)\cap A^{p,q}(X)
\end{align*}

We have

\begin{theorem}  (Hodge decomposition theorem) Let $X$ be a 
compact K\"ahler manifold of dimension $n$. Then the n-th 
complex de Rham cohomology group of $X$ can be written as the
direct sum
\begin{equation*}
H^n(X,Z)\otimes { C}
=H^n_{DR}(X,C)=
\underset{p+q=n}{\oplus}H^{p,q}(X)
\end{equation*}
\end{theorem}
\vspace{0.2in}

Suppose $\omega$ is the K\"ahler form of $X$.
Define
\[
L: H^k(X,C)\rightarrow H^{k+2}(X,C),
\quad
k=0,\cdots, 2n-2
\]
to be the multiplication by $\omega$.

Suppose $X$ is a smooth algebraic 
variety. Then $\omega$ is called a polarization of $X$ 
if $[\omega]$ is the first Chern class of an ample line bundle over $X$.
In that case, 
 $(X,\omega)$ is
called the  polarized algebraic variety.

The following two famous Lefschetz theorems are extremely important
in defining the classifying space and the period map.

\begin{theorem} (Hard Lefschetz theorem) On a polarized algebraic variety
$(X,\omega)$ of dimension $n$,
\[
L^k: H^{n-k}(X,C)\rightarrow H^{n+k}(X,C)
\]
is an isomorphism for every positive integer $k\leq n$.
\end{theorem}

The primitive cohomology $P^k(X,C)$ is  then defined to be
the kernel of $L^{n-k+1}$ on $H^k(X,C)$.

\begin{theorem}(Lefschetz Decomposition Theorem)
On a polarized algebraic variety $(X,\omega)$, we have 
the following decomposition:
\[
H^n(X,C)=\oplus_{k=0}^{[\frac n2]} L^kP^{n-2k}(X,C)
\]
\end{theorem}

\vspace{0.2in}

Define $H_Z=P^n(X,C)\cap H^n(X,Z)$ and
$H^{p,q}=P^n(X,C)\cap H^{p,q}(X)$. Then we have
\[
H_Z\otimes C=\sum H^{p,q}, \quad H^{p,q}=\bar{H^{q,p}}
\]
for $p+q=n$. Set $H=H_Z\otimes C$. We call $\{H^{p,q}\}$ the Hodge
decomposition of $H$.

\begin{rem}
We  define a filtration of $H_Z\otimes C=H$ by
\[
0\subset F^n\subset F^{n-1}\subset \cdots F^1=H
\]
such that
\[
H^{p,q}=F^p\cap\bar{F}^q,\qquad F^p\oplus\bar{F^{n-p+1}}=H
\]
The set $\{H^{p,q}\}$ and $\{F^p\}$ are equivalent to describe
the Hodge decomposition of $H$. In the remaining of this paper, we 
will use both notations interchangeably.
\end{rem}

Now suppose that $Q$ is the quadratic form on $H_Z$ 
induced by the cup product of the cohomology group.  
$Q$ can be represented by
\[
Q(\phi, \psi)=(-1)^{n(n-1)/2}\int \phi\wedge \psi
\]
for $\phi, \psi\in H$.
$Q$ is a nondegenerated form, and is skewsymmetric if $n$ is odd 
and is symmetric if $n$ is even. It satisfies the two 
Hodge-Riemannian relations
\begin{enumerate}\label{hr}
\item $Q(H^{p,q},H^{p',q'})=0 \quad unless\quad p'=n-p, q'=n-q$;
\item $(\sqrt{-1})^{p-q}\,Q(\phi,\bar\phi)>0$ for any nonzero 
element $\phi\in H^{p,q}$.
\end{enumerate}

\begin{definition}
A polarized Hodge structure
of weight $n$, denoted by \linebreak
$\{H_Z,F^p,Q\}$, is given by a 
filtration of $H=H_Z\otimes C$
\[
0\subset F^n\subset F^{n-1}\subset\cdots\subset F^0\subset H
\]
such that
\[
H=F^p\oplus\bar{F}^{n-p+1}
\]
together with a bilinear form
\[
Q: H_Z\otimes H_Z\rightarrow Z
\]
which is skewsymmetric if $n$ is odd and symmetric if $n$ is 
even such that it satisfies the two Hodge-Riemannian 
relations:

3.  $Q(F^p,F^{n-p+1})=0$ for $p=1,\cdots n$;

4.  $(\sqrt{-1})^{p-q}\,Q(\phi,\bar\phi)>0$ if $\phi\in 
H^{p,q}$ and $\phi\neq 0$.

where $H^{p,q}$ is defined by
\[
H^{p,q}=F^p\cap\bar{F}^q
\]
for $p+q=n$.
\end{definition}

\begin{definition}
The classifying space $D$ for the 
polarized Hodge structure is the set of all the filtration
\[
0\subset F^n\subset\cdots\subset F^1\subset H, \qquad
F^p\oplus\bar{F^{n-p+1}}=H
\]
or the set of all the decompositions
\[
\sum H^{p,q}=H,\qquad H^{p,q}=\bar{H^{q,p}}
\]
on which $Q$ satisfies the two
Hodge-Riemannian relations 1,2 or 3,4 above.
\end{definition}

Over the classifying space $D$ we have the holomorphic vector bundles
$\underline{F}^n,
\cdots, \underline{F}^1, \underline{H}$ whose fibers at each point are the
vector spaces $F^n,\cdots, F^1,H$, respectively. These bundles are called
Hodge bundles.

It is well known that the holomorphic tangent bundle $T(D)$ can
be realized as a subbundle of $Hom(\underline H,\underline H)$:
\[
T(D)\subset \oplus Hom(\underline F^p,\underline H/\underline F^p)=
\underset{r> 0}{\oplus}
Hom( \underline H^{p,q},\underline H^{p-r,q+r})
\]
such that the following compatible condition holds

\[
\begin{CD}
F^{p} @> >>F^{p-1}\\
@VVV        @VVV\\
H/F^p@ <<< H/F^{p-1}
\end{CD}
\]

\begin{definition}  
A subbundle $T_h(D)$  called the
horizontal bundle of $D$, if
\[
T_h(D)=\{\xi\in T(D)| \xi F^p\subset F^{p-1}, p=
1,\cdots, n\}
\]
\end{definition}

For any point $p\in D$ such that $p$ is defined as subspaces 
$\{H^{p,q}\}$ of $H$, define the two vector spaces
\begin{align*}
& H^+=H^{n,0}+H^{n-2,2}+\cdots\\
& H^-=H^{n-1,1}+H^{n-3,3}+\cdots
\end{align*}
Now we fix a point $p_0\in D$. Suppose the corresponding vector spaces are 
$\{H_0^{p,q}\}$ and $\{H_0^+, H_0^-\}$. We define $V$ to be the connected
compact subgroup of $G$ leaving $\{H_0^{p,q}\}$ unchanged for any
$p,q$ with $p+q=n$, and $K$ to be the connected compact subgroup of $G$
leaving $H_0^+$
invariant. We
give the basic properties of the classifying spaces in the following
three lemmas. 
The proofs are easy and are omitted.

\begin{lemma}
\label{121}
$K$ is the maximal compact subgroup of 
$G$ containing $V$. Thus $V$ itself is a compact subgroup.
\end{lemma}
\qed

We define the Weil operator
\[
C: H^{p,q}\rightarrow H^{p,q},\quad C|_{H^{p,q}}=(\sqrt{-1})^{p-q}
\]
Then we have
\[
C|_{H^+}=(\sqrt{-1})^n,\qquad C|_{H^-}=-(\sqrt{-1})^n
\]
 
Let
\[
Q_1(x,{y})=Q(Cx,\bar{y})
\]
 
Then we have

\begin{lemma}
$Q_1$ is an Hermitian inner product.
\end{lemma} \qed
\begin{lemma}
Let
\[
D_1=\{H^{n,0}+H^{n-2,2}+\cdots|\{H^{p,q}\}\in D\}
\]
Then the group $G$ acts on $D_1$ transitively with the stable subgroup 
$K$ at $H_0^+$, and  $D_1$ is a symmetric space.
\end{lemma}
\qed

\begin{definition}
We call map $p$
\[
p: G/V\rightarrow G/K,\qquad \{H^{p,q}\}\mapsto H^{n,0}+H^{n-2,2}+\cdots
\]
the natural projection of the classifying space.
It is the same as the map $aV\mapsto aK$ between the coset spaces of $V$
and $K$.
\end{definition}

With the above discussions, we can prove

\begin{prop}\label{pp21}
Suppose $T_v(D)$ is the distribution of the tangent vectors of the fibers 
of the canonical map
\[
p:\, D\rightarrow G/K
\]
then
\[
T_v(D)\cap T_h(D)=\{ 0\}
\]
\end{prop}

{\bf Proof:}
Let ${\frk g}$ be the Lie algebra of the Lie group $G$. Let ${\frk
g}={\frk f}+{\frk p}$ be the Cartan decomposition such that ${\frk f}$ is
the Lie algebra of $K$. Then
\[
T_v(D)=G\times_V{\frk v}_1
\]
where ${\frk f}={\frk v}+{\frk v}_1$ and ${\frk v}_1$ is the orthonormal 
complement of the Lie algebra ${\frk v}$ of $V$. On the other
hand, $T_h(D)\subset G\times_V{\frk p}$. So we have 
$T_v(D)\cap T_h(D)=\{ 0\}$.

 \qed

\begin{definition}
A horizontal slice ${\mathcal M}$ of $D$ is a complex integral submanifold
of the distribution $T_h(D)$.
\end{definition}

\begin{definition}
Let $U$ be an open neighborhood of the universal deformation space at a
Calabi-Yau manifold $X$. Then for each $X'$ near $X$, we have an
isomorphism $H^n(X',C)=H^n(X,C)$. Under this isomorphism, 
$\{H^{p,q}(X')\cap P^n(X',C)\}_{p+q=n}$ can be considered as a point of
$D$. The map 
\[
U\rightarrow D\qquad X'\mapsto \{H^{p,q}(X')\cap P^n(X',C)\}_{p+q=n}
\]
is called the period map.
\end{definition}

The most important property of the period map is the following~\cite{Gr}:

\begin{theorem} (Griffiths)
The period map $p: U\rightarrow D$ is an immersion and  
$p(U)$ is a horizontal slice of the classifying space.
\end{theorem}

From the above theorem and Proposition~\ref{pp21} in this section, we can
prove:

\begin{cor}\label{cor31}
With the notations as above, the map
\[
p: U\subset D\rightarrow G/K
\]
is an immersion.
\end{cor}

\section{The Invariant
Complex Structure}

It was proved by Griffiths that 
the classifying space $D$ is actually  a homogeneous complex
manifold. We are going to study this fact a little bit in detail
via the Lie group point of view.

Define $H_R=H_Z\otimes R$ and let
\begin{align*}
G_R&=Aut(H_R,Q)
=\{g:H_R\rightarrow H_R|Q(g\phi,g \psi)=g(\phi,\psi),
\phi,\psi\in H_R\}\\
G_C&=Aut(H_C,Q)
=\{g:H_R\rightarrow H_R|Q(g\phi,g \psi)=g(\phi,\psi),
\phi,\psi\in H\}
\end{align*}

Let $G=G_R$.
Then $G$ acting on $D$ transitively and thus $D$ is a
homogeneous space.

Let $V$ be the isotropy group fixing one point of $D$,
 then $V$ is a compact subgroup (see Lemma~\ref{121}).

\smallskip

Let ${\frk g}$ be the  Lie algebra
of $G$, then we have the standard Cartan decomposition
\[
{\frk g}={\frk f}+{\frk p}
\]
into the compact part ${\frk f}$ and noncompact part
 ${\frk p}$.
We assume that the Lie algebra
${\frk v}$ of $V$ is contained in
${\frk f}$.
 Since the Killing form on ${\frk f}$ is
negative definite, there is a subset
${\frk v}_1$ of ${\frk f}$ such that
\[
{\frk f}={\frk v}+{\frk v}_1
\]
is an orthnormal decomposition of ${\frk f}$. Thus we have
\[
{\frk g}={\frk v}\oplus {\frk v}_1\oplus {\frk p}
\]
There is a natural representation of ${\frk v}$ to ${\frk v}_1\oplus
{\frk p}$ such that the tangent bundle of $D$ is the associated
bundle of the principle bundle $G\rightarrow G/V$ with respect
to this representation.

Suppose that ${\frk a}^c$ is the complexification of
 a Lie algebra ${\frk a}$. Then we have
  \[
{\frk g}^c={\frk v}^c\oplus {\frk v}_1^c\oplus {\frk p}^c
\]

Let $\tau$ be the complex conjugate of ${\frk g}^c$ with respect to the
compact real form of ${\frk g}^c$. Then
\[
{\frk v}_1^c\oplus {\frk p}^c={\frk n}_-\oplus\tau({\frk n}_-)
\]
be the splitting into $(\pm \sqrt{-1})$ $\tau$-spaces.
Then complex structure of $D$ is determined by $\tau({\frk n}_-)$.

\smallskip

Suppose $J$ is the invariant complex structure of $D$. It is 
well known that to give an invariant complex structure on $G/V$  
is equivalent to  give a
 linear transformation $J$ on ${\frk v}_1+{\frk 
p}$ such that
\begin{align*}
& J^2=-id_{({\frk v}_1+{\frk p})}\\
& \rho(h)J=J\rho(h),\qquad\forall h\in V
\end{align*}
where
\[
\rho : V\rightarrow {\frk v}_1+{\frk p}
\]
is the standard adjoint representation.

By the structure theory of the complex 
semi-simple Lie algebra that the complexification of ${\frk 
v}_1+{\frk p}$ can be written as the sum of the root spaces
\[
({\frk v}_1+{\frk p})^c=\sum_{\alpha\in I}{\frk g}_{\alpha}
\]
for some index set $I$. Suppose ${\frk h}$ is the Cartan 
subalgebra of ${\frk g}^c$. Let ${\frk h}_R$ be its real part, then
\[
{\frk h}_R\subset {\frk v}
\]
thus $J$ is $V$-invariant implies that $J$ is ${\frk h}$ 
invariant. In particular
\[
[h,JX]=J[h,X]
\]
for $h\in{\frk h}, X\in ({\frk v}_1+{\frk p})^c$. Now let $X\in 
{\frk g}_\alpha$, then 
\[
[h,JX]=J[h,X]=J\alpha (h)X=\alpha (h) JX
\]
thus
\[
JX\in{\frk g}_{\alpha}
\]

Since $J^2=-1$, we know $JX=\pm\sqrt{-1}X$. In particular, if 
$J_1$, $J_2$ are the two invariant complex structures, we have
\[
J_1J_2=J_2J_1
\]

Suppose $p$ is the projection of ${\frk v}_1+{\frk p}$ to the 
second factor. 

\begin{definition}\label{df41}
We call $\omega$ is the fundamental form on $D$ if $\omega$ is 
$G$ invariant and if on ${\frk v}_1+{\frk p}$, the tangent space of $D$
at $T_e(D)$,  $\omega$ is 
defined as 
\[
\omega (X,Y) = - B(pX,pJY)
\]
where $B$ is the Killing form of the Lie algebra.
\end{definition}

From the following proposition, we know
$\omega$ is well defined.
\begin{prop}
We have
\begin{enumerate}
\item $\omega (X,Y) =-\omega (Y,X), for \quad X,Y\in
{\frk v}_1+{\frk p}$;
\item $\omega$ is $V$ invariant.
\end{enumerate}
\end{prop}

{\bf Proof:}
We recall the root decomposition
\[
({\frk v}_1+{\frk p}_0)^c=\sum_{\alpha\in I}{\frk g}_{\alpha}
\]
and
\[
J{\frk g}_\alpha\subset{\frk g}_\alpha\quad for \quad\alpha\in I
\]
In particular, if $X\in {\frk v}_1, JX\in{\frk v}_1$; and  $X\in{\frk
p}, 
JX\in{\frk p}$.

So if $X\in{\frk v}_1$, then $pX=0$, $pJX=0$. So
\[
\omega (X,Y)=0=-\omega (Y,X)
\]

In order to verify the first claim, 
 we need only to check the proposition for $X\in {\frk 
g}_{\alpha}$ and $Y\in{\frk g}_{\beta}$ for noncompact roots 
$\alpha$ and $\beta$.
Suppose
$JX=\sigma_XX$, 
$JY=\sigma_YY$.
Then
$\omega(X,Y)=\sigma_YB(X,Y)$,
$\omega(Y,X)=\sigma_XB(Y,X)$.
If $\alpha+\beta\neq 0$, then $B(X,Y)=0$.
So we need only assume that $\alpha+\beta=0$. In this case, 
suppose that
$JX=\sigma X$. Then
$J\bar{X}=-\sigma\bar{X}$.
But $\bar{X}\in{\frk g}_{-\alpha}$.
So $JY=-\sigma Y$
Thus  we have
$\sigma_X=-\sigma_Y$
and then
$\omega(X,Y)=-\omega(Y,X)$.

Next, $\forall h\in V$, we have
\[
\omega(Ad(h)X,Ad(h)Y)=\omega(X,Y)
\]
because $Ad(h)$ commutes $p$ and $J$.

From the above theorem we know that $\omega$ is well defined.

\section{The Pluriharmonicity}
In this section, we are going to prove Theorem~\ref{thm1}.
The notations are as in the previous sections.

Let $D=G/V$ be a classifying space and ${\frk g}$ be the Lie algebra of
$G$. We begin with the following key observation:

\begin{theorem}
\label{thm311}
Let ${\frk g}={\frk v}+{\frk v}_1+{\frk p}$ be the  
decomposition in the
previous sections. 
Suppose $U$ is an open neighborhood
of a horizontal slice.
 Then 
 if $X,Y,Z\in\Gamma(U,G\times_V{\frk p})$ or
$X,Y\in\Gamma(U,G\times_V{\frk v}_1)$ and $Z\in\Gamma(U,G\times_V{\frk
p})$, then
\[
d\omega(X,Y,Z)=0
\]
where $\omega$ is the differential form defined in Definition~\ref{df41}.
\end{theorem}

{\bf Proof:}
We have
\[
T(G/V)=T_v(G/V)+G\times_V{\frk p}
\]
where
\[
T_v(G/V)=G\times_V{\frk v}_1
\]

Suppose $\sigma_0 : G\rightarrow G$ is the involution. 
That is , $\sigma_0$ is a isomorphism of $G$ such that $\sigma_0^2=1$. 
$\sigma_0$ induced the Lie algebra isomorphism
$\sigma_0: {\frk g}\rightarrow {\frk g}$. It is easy to see that
 $\sigma_0(X)=X$ for $X\in{\frk f}$,
$\sigma_0(X)=-X$ for $X\in{\frk p}$.

$\sigma_0$  also induced a $C^\infty$ map
\[
G/V\rightarrow G/V, \qquad aV\mapsto \sigma(a)V
\]
Let $\sigma_g$ be the map 
\begin{equation}\label{sym}
\sigma_g=(L_g)\sigma_0(L_{g^{-1}})
\end{equation}
on $G/V$ where $L_g$ is the left transformation, then
$\sigma_g(gV)=gV$.

\begin{lemma}\label{lm41}
$\sigma_{g}(\omega)=\omega$ where $\omega$ is the differential 
form defined in Definition~\ref{df41}.
\end{lemma}

{\bf Proof:} By the definition of $\sigma_g$, we need only 
prove that $\sigma_0(\omega)=\omega$.
Next we observe that $\sigma_0(\omega)$ is also an invariant 
form. So we only check that $\sigma_0(\omega)=\omega$
at the original point.

We see that
\[
\sigma_0J=J\sigma_0, p\sigma_0=\sigma_0p
\]
where $p:D\rightarrow D_1$ is the projection.

So for any tangent vector $X,Y$,
\begin{align*}
&(\sigma_0\omega)(X,Y)=\omega(\sigma_0(X),\sigma_0(Y))
=-B(p\sigma_0X,pJ\sigma_0Y)\\
&=-B(\sigma_0pX,\sigma_0pJY)=-B(pX,pJY)
=\omega(X,Y)
\end{align*}
where $B$ is the Killing form. Thus the lemma is proved.

\qed

{\bf Continuation of the Proof of Theorem ~\ref{thm311}:} 

Now if $X,Y,Z\in\Gamma(U,G\times_V{\frk p})$ or
$X,Y\in\Gamma(U,G\times_V{\frk v}_1)$ and $Z\in\Gamma(U,G\times_V{\frk
p})$, then

\begin{align*}
& d\omega(X,Y,Z)=-d\omega(\sigma_0 X,\sigma_0 Y,\sigma_0 Z)\\
&=-(\sigma_0d\omega)(X,Y,Z)=-d(\sigma_0\omega)(X,Y,Z)\\
&=-d\omega(X,Y,Z)
\end{align*}

So
\[
d\omega(X,Y,Z)=0
\]

\qed

The invariant Riemannian metric on $G/K$ 
is defined by (up to a multiplication of a constant),
\[
(X,Y)_{gV}=B((L_{g^{-1}})_*X, (L_{g^{-1}})_*Y)
\]
where $L_{g^{-1}}$ is the left translation of the group $G$ by
$g^{-1}$.

\begin{lemma}
\label{thm21}
Let $X,Y$ be vector fields on $D$. Let $p: D\rightarrow D_1$ be the
projection. Then
\[
\omega (X,Y)=-(p_*X,p_*JY)
\]
\end{lemma}

{\bf Proof:}
At the origin point, the lemma is trivially true by 
definition. At a general point $gV$ of $D$, note that the left 
translation $L_{g^{-1}}$ commutes with the projection $p$, we have
\begin{align*}
&\omega_{gV}(X,Y)=\omega_0((L_{g^{-1}})_*X,(L_{g^{-1}})_*Y)
=-B(p(L_{g^{-1}})_*X,pJ(L_{g^{-1}})_*Y)\\
&=-B((L_{g^{-1}})_*p_*X,(L_{g^{-1}})_*p_*JY)
=-(p_*X,p_*JY)
\end{align*}

\qed

\begin{definition}
An immersion 
\[
p: M\rightarrow N
\]
from a complex manifold $M$ to a Riemannian manifold $N$
is called pluriharmonic, if for any vector field $X$ on $M$
\[
\nabla_{p_*X}p_*X+\nabla_{p_*JX}p_*JX+p_*J[X,JX]=0
\]
for the Riemannian connection $\nabla$ of $N$.
\end{definition}

Now we begin to prove Theorem~\ref{thm1}:

\begin{theorem}
\label{thmfud}
Suppose ${U}$ is a horizontal slice of $D$.  
Then
\[
p: {U}\subset D\rightarrow G/K
\]
is a pluriharmonic immersion for the invariant Riemannian connection
$\nabla$ on
$G/K$.
\end{theorem}

{\bf Proof:} 
We have already proved that the map is an immersion 
in the Corollary~\ref{cor31}. So it remains to prove the
pluriharmonicity of the
map.
We say a complex submanifold $S$ is integrable {\it at a point $q$}, if 
\[
T_q(S)\subset (G\times_V{\frk p})_q
\]
and if at a neighborhood of $q$, the map $p$ is an immersion
at $q$ from $S$ to $G/K$.
Let  $U$ be    an open neighborhood of $S$. Let
$X,Y,Z\in\Gamma(U,T{U})$ and $X_q,Y_q,Z_q\in (G\times_V{\frk
p})_q$. Then by Theorem~\ref{thm311} at $q$,
\[
d\omega(X,Y,Z)=0
\]
$p_*X,p_*Y,p_*Z$ are well 
defined and 
$C^{\infty}$ in a neighborhood of $U$. Then by Lemma ~\ref{thm21}
\begin{align*}
0&=-d\omega(X,Y,Z)=-X\omega(Y,Z)+Y\omega(X,Z)-Z\omega(X,Y)\\
&-\omega(X,[Y,Z])+\omega(Y,[X,Z])-\omega(Z,[X,Y])\\
&=(p_*X)(p_*Y,p_*JZ)-(p_*Y)(p_*X,p_*JZ)+(p_*Z)(p_*X,p_*JY)\\
&-(p_*JX,[p_*Y,p_*Z])+(p_*JY,[p_*X,p_*Z])-(p_*JZ,[p_*X,p_*Y])\\
&=(\nabla_{p_*X}p_*Y,p_*JZ)+(p_*Y,\nabla_{p_*X}p_*JZ)
-(\nabla_{p_*Y}p_*X,p_*JZ)\\
&-(p_*X,\nabla_{p_*Y}p_*JZ)+(\nabla_{p_*Z}p_*X,p_*JY)
+(p_*X,\nabla_{p_*Z}p_*JY)\\
&-(p_*JX,\nabla_{p_*Y}p_*Z)+(p_*JX,\nabla_{p_*Z}p_*Y)
+(p_*JY,\nabla_{p_*X}p_*Z)\\
&-(p_*JY,\nabla_{p_*Z}p_*X)-(p_*JZ,\nabla_{p_*X}p_*Y)
+(p_*JZ,\nabla_{p_*Y}p_*X)\\
&=(p_*Y,\nabla_{p_*X}p_*JZ)+(p_*JY,\nabla_{p_*X}p_*Z)
-(p_*X,\nabla_{p_*Y}p_*JZ)\\
&+(p_*X,\nabla_{p_*Z}p_*JY)-(p_*JX,\nabla_{p_*Y}p_*Z)
+(p_*JX,\nabla_{p_*Z}p_*Y)
\end{align*}

If we substitute $X$ by $JX$, $Y$ by $JY$ and $Z$ by $JZ$, we have
\begin{align*}
0&=-(p_*JY,\nabla_{p_*JX}p_*Z)-(p_*Y,\nabla_{p_*JX}p_*JZ)
+(p_*JX,\nabla_{p_*JY}p_*Z)\\
&-(p_*JX,\nabla_{p_*JZ}p_*Y)+(p_*X,\nabla_{p_*JY}p_*JZ)
-(p_*X,\nabla_{p_*JZ}p_*JY)
\end{align*}

Substituting $X$ by $JX$, we have
\begin{align*}
0&=(p_*Y,\nabla_{p_*JX}p_*JZ)+(p_*JY,\nabla_{p_*JX}p_*Z)
-(p_*JX,\nabla_{p_*Y}p_*JZ)\\
&+(p_*JX,\nabla_{p_*Z}p_*JY)+(p_*X,\nabla_{p_*Y}p_*Z)
-(p_*X,\nabla_{p_*Z}p_*Y)
\end{align*}

Comparing the above two equations, we get
\begin{align*}
0&=-2(p_*JY,\nabla_{p_*JX}p_*Z)-2(p_*Y,\nabla_{p_*JX}p_*JZ)
+(p_*JX,[p_*JY,p_*Z])\\
&-(p_*JX,[p_*JZ,p_*Y])
+(p_*X,[p_*JY,p_*JZ])-(p_*X,[p_*Y,p_*Z])
\end{align*}
where the last four terms is equal to
\begin{align*}
&-(p_*X,p_*J[JY,Z])+(p_*X,p_*J[JZ,Y])\\
&+(p_*X,p_*[JY,JZ])-(p_*X,p_*[Y,Z])=0
\end{align*}
because of the integrability condition of $J$.

So we have
\begin{equation}
\label{eq1}
(p_*Y,\nabla_{p_*JX}p_*JZ)+(p_*JY,\nabla_{p_*JX}p_*Z)=0
\end{equation}
Let $X=Z$, we have
\[
(p_*Y,\nabla_{p_*JX}p_*JX)+(p_*JY,\nabla_{p_*JX}p_*X)=0
\]

Substituting $X$ by $JX$, we have
\[
(p_*Y,\nabla_{p_*X}p_*X)-(p_*JY,\nabla_{p_*X}p_*JX)=0
\]

So
\[
(p_*Y,\nabla_{p_*X}p_*X)+(p_*Y,\nabla_{p_*JX}p_*JX)
+(p_*JY,p_*[JX,X])=0
\]
Then the theorem follows from the fact that for any $\tilde X,\tilde
Y,\tilde Z\in T_{p(q)}(G/K)$, there is an integral submanifold $S$ at $q$
and $X,Y,Z\in T_q(D)$ such that $p_*X=\tilde X$, $p_*Y=\tilde Y$,
$p_*Z=\tilde Z$.

\begin{cor}
Use the same notation as above,
We also have
\[
(p_*Z,\nabla_{p_*JX}p_*JY)
+(p_*JZ,\nabla_{p_*JX}p_*Y)=0
\]
\end{cor}

{\bf Proof:} This is the same as Equation~\ref{eq1}.

\section{The Hodge Metric}
We use the notations in the previous 
sections. Now suppose $U$ 
is a horizontal slice. 
 We also denote $\omega$ the restriction 
to $U$ of the 
fundamental form in Definition~\ref{df41}. 
Then we have

\begin{lemma}
On the horizontal slice $U$
\begin{enumerate}
\item  $d\omega=0$;
\item $\omega$\ is a 1-1 form;
\item $\omega(X,JX)>0, X\ne 0$
\end{enumerate}
\end{lemma}

{\bf Proof:} The first assertion is easy from 
Theorem ~\ref{thm311}.
 To prove the second assertion, let 
assume
$X_\alpha, X_\beta\in T^{1,0}(G/V)$, where $\alpha,\beta$ 
are roots of ${\frk g}^c$. Then $\omega(X_\alpha,X_\beta)\neq 0$
if and only if $\beta=-\alpha$. In that case, we can assume
\[
X_\beta=-\bar{X}_\alpha
\]
contradicting to the fact that both of them are in $T^{1,0}(G/V)$. So
$\omega$ is a 1-1 form on $G/V$, thus a 1-1 form 
on $U$. If $X\neq 0$,
\[
\omega(X,JX)=(p_*X,p_*X)>0
\]
if $X\neq 0$ by Lemma~\ref{thm21}.

Thus $\omega$ defined a K\"ahler metric whose underlying 
Riemannian metric is the restriction of the invariant Riemannian metric on
$U$.

\begin{definition}
The K\"ahler metric $\omega$ is called the Hodge metric of the horizontal
slice $U$.
\end{definition}

The main theorem of this section is

\begin{theorem}
With respect to the Hodge metric, the holomorphic bisectional curvature
and
the Ricci curvature are nonpositive. Furthermore, the holomorphic sectional curvature and
the Ricci curvature are negative away from zero by a constant number.
\end{theorem}

{\bf Proof:}
By Lemma~\ref{thm21}, the Hodge metric is the restriction of the invariant
Riemannian metric on $G/K$.
Suppose $R,\tilde R$ are the curvature tensors of $G/K$ and $U$, 
respectively. Then we have the Gauss formula:
\begin{align*}
& \tilde R(X,Y,X,Y)+\tilde R(X,JY,X,JY)\\
&=R(X,Y,X,Y)+R(X,JY,X,JY)\\
&+(h(X,X),h(Y,Y)+h(JY,JY))
-|h(X,Y)|^2-|h(X,JY)|^2
\end{align*}
where $h(\cdot,\cdot)$ refers to the second fundamental form and $X,Y\in
TU$. Now
the pluriharmonicity of the map 
(Theorem~\ref{thm1}) $p: U\rightarrow G/K$ implies that
\[
h(Y,Y)+h(JY,JY)=0
\]
Thus the negativity of the bisectional curvature follows from the
following lemma:

\begin{lemma}
\label{lem123}
$G/K$ has nonpositive sectional curvature. Furthermore, at the original
point
\[
R(X,Y,X,Y)=-||[X,Y]||^2
\]
\end{lemma}

{\bf Proof:} 
Since $G/K$ is a symmetric space of noncompact type,
the curvature formula follows from~\cite{Hel}. in particular, the curvature tensor is
nonpositive.

\qed

{\bf Continuation of the Proof of the Theorem:} 
From the above lemma, we see that the holomorphic bisectional curvature is
nonpositive. 
Now we consider the
holomorphic sectional
curvature. By the above computation, we have
\[
\tilde R(X,JX,X,JX)\leq R(X,JX,X,JX)
\]

That $R(X,JX,X,JX)$ is negative away from zero is hinted by Griffiths and
Schmidt~\cite{GS}: 
if we can prove that for any $X\in T_eU$, 
\[
|[X,JX]|\geq c||X||^2
\]
for $c>0$, then we have
\[
R(X,JX,X,JX)\leq -c^2||X||^4
\]

Thus the holomorphic sectional curvature is negative away from zero at the original
point. Because the homegeneity of the metric, we know that it is negative everywhere.

Now we let
\[
X=\sum_i a_ig_i+\bar{a_i}g_{-i}
\]
where $Jg_i=\sqrt{-1}g_i$ and $g_{-j}=\bar{g_j}$. Then
\[
[X,JX]=-2\sqrt{-1}\sum_i |a_i|^2[g_i,g_{-i}]+\cdots
\]
here $\cdot$ refers to the terms which are not in the maximum torus. 
Since $\sum_i
|a_i|^2[g_i,g_{-i}]$ belongs to the cone of positive roots, 
it will be never zero unless
$a_i\equiv 0$. Thus we proved that $||[X,JX]||\geq c||X||^2$. 

The Ricci curvature is the sum of the holomorphic bisectional curvature of
different directions. So its negativity follows from the negativity of the
holomorphic sectional curvature.

One application of the above theorem is that we might be able to use it to
give a proof the compactification theorem of the moduli space. It is easy
to see that the Hodge metric can be defined on the moduli space as a
K\"ahler orbifold metric. With this, we can proved that

\begin{cor}
If ${\mathcal M}$ is complete and smooth  with bounded sectional
curvature with respect to the Hodge metric, then ${\mathcal M}$
is quasi-projective.
\end{cor}

{\bf Proof:} It follows from a theorem of Yeung~\cite{Yeung1}.

\begin{rem}
The above corollary is valid in a much weaker assumption. We will write out the proof in
a subsequent paper.
\end{rem}

\section{On $\Gamma$-invariant K\"ahler Metrics}
Let $D=G/V$ and let $\Gamma$ be a discrete subgroup of $G$. We say that
$\Gamma$ is
cocompact if $\Gamma\backslash G$ is a compact topological space.
In this section, we restrict ourselves to the case that 
$D$ is the classifying space of some Calabi-Yau
threefold.
We consider the problem of the existence of K\"ahler metrics on $D$ which
are
$\Gamma$-invariant.

The motivation of this problem is from the recent work of
Rajan~\cite{Raj}. In his paper, Rajan proved the infinitesmal rigidity of
the complex structure for a large family of homogeneous K\"ahler
manifolds. It would be very interesting to prove the same rigidity theorem
for the classifying space. However, In general, $D$ is not a homogeneous
K\"ahler manifold. Furthermore, we proved 
in this section that $D$ even does
not admit a K\"ahler metric which is $\Gamma$-invariant.

To be more precise, we have

\begin{theorem} Suppose $D$ is the classifying space
of a Calabi-Yau threefold. $D$ is the set of filtration
\[
0\subset F^3\subset F^2\subset F^1\subset H
\]
satisfying the two  Hodge-Riemannian conditions.
Suppose  that  $\Gamma$ is a cocompact lattice of $D$. Then  there
are no 
K\"ahler metrics on $D$ which are $\Gamma$-invariant.
\end{theorem}

Recall that $\Gamma$ is cocompact if and only if $\Gamma\backslash D$
is a compact space.
The idea to prove the theorem is to  consider the 
map $p$ of the projection
\[
p: D\rightarrow  G/K
\]
where $K$ is the maximum compact subgroup of $G$ containing $V$.

It is easy to see that
$G={\frk S}{\frk p}(2n+2,R)$,
the real symplectic group,
where $2n+2$ is the complex dimension of $H$. 
 In particular, $G/K$ is the Hermitian symmetric
space, which can be realized as the  Siegal Space of the third kind.

We observed that $p$ is not holomorphic nor anti-holomorphic. The reason is
that if $p$ is anti-holomorphic, then we can reverse the complex structure
on $D$ to make it holomorphic. By a general theorem in
Murakami~\cite{M}, we know that in that case, $D$ will be a
homogeneous K\"ahler manifold. But we have

\begin{lemma}
 $D$ is not a
homogeneous K\"ahler manifold.
\end{lemma}

{\bf Proof:} The proof is easy after written out the root decomposition
and the conjugate $\tau$ explicitly. For details, see~\cite{Lu}.

\qed

\begin{lemma}
If $D$ admits a $\Gamma$-invariant K\"ahler metric, then there
is a $\Gamma$-equivariant holomorphic or anti-holomorphic  map $f:
D\rightarrow G/K$
which is surjective at some point $x_0\in D$.
\end{lemma}

{\bf Proof:} Suppose $\tilde\omega$ is the K\"ahler form on $D$ which is
$\Gamma$-invariant.
According to Eells-Sampson~\cite{ES}, Jost-Yau~\cite{JY2}, and
Labourie~\cite{L}, 
there is a $\Gamma$-equivariant harmonic map
$f:D\rightarrow G/K$ such that
$f$  is
$\Gamma$-equivariant
 homotopic to $p$. 
That is, $f$ and $p$ can be linked by a path of continuous
$\Gamma$-equivariant maps. 
By topology, we know $p,f$ induce the same 
homomorphism between the cohomology groups:
\[
p^*, f^*: H^{2N}(\Gamma\backslash G/K,C)\rightarrow
H^{2N_1}(\Gamma\backslash D, C)
\]

Here  $2N$ is the real dimension of $\Gamma\backslash G/K$. $2N_1$ is
the real dimension of $\Gamma\backslash D$.
In particular,  if  $\eta$ is the volume
form of $\Gamma\backslash G/K$ then we have

\[
\int_{\Gamma\backslash D}f^*\eta\wedge\tilde\omega^{N_1-N}
=\int_{\Gamma\backslash D}p^*\eta\wedge\tilde\omega^{N_1-N} 
\]

On the other hand
\[
\int_{\Gamma\backslash D}p^*\eta\wedge\tilde\omega^{N_1-N}\
=\int_{\Gamma\backslash G/K}\eta\int\tilde\omega^{N_1-N}>0
\]

So 
\[
\int_{\Gamma\backslash D}f^*\eta\wedge\tilde\omega^{N_1-N}\neq 0
\]
and in particular, $f^*\eta\neq 0$ at some point $x_0$.

By the rigidity theorem of Siu~\cite{Siu1}, $f$ is a holomorphic or anti-holomorphic  
map.

\qed

The invariant Riemannian metric $g$ on $D$ is a Hermitian metric.

\begin{lemma}
The map $f$ is a harmonic map with respect to $g$.
\end{lemma}

{\bf Proof:} We assume that $f$ is holomorphic without loosing generality.
 First note that each fiber of $p: D\rightarrow G/K$ is a compact K\"ahler
submanifold 
of $D$. Thus according to Liouville's theorem $f$ is a constant along 
each fiber. 

Now suppose that $(w^1,\cdots,w^N)$ is the holomorphic local coordinate of
$G/K$ and $(z^1,\cdots
z^{N_1})$ is the corresponding holomorphic local coordinate such that 
$f=(f^1,\cdots, f^N)$ with respect to these coordinates are holomorphic
functions. 
Thus at any  point $x\in D$,
\[
\frac{\pa f^\alpha}{\pa\bar{z}^j}(x)=0
\]
for any $\alpha, j$. 
We assume that at $p(x)$, 
$(w^1,\cdots, w^N)$ is normal. That is, all the connection coefficients
are zero at $p(x)$.

Now suppose $\Delta$ is the Laplacian of the
metric $g$.
We need to prove
\[
\Delta f^\alpha(x)=0
\]
for $\alpha=1,\cdots, N$.

Now 
\begin{align*}
&-\Delta f^\alpha=\delta d f^\alpha=\delta \pa f^\alpha
+\delta\bar\pa f^\alpha\\
&=-*d*\frac{\pa f^\alpha}{\pa z^i} dz^i
=-*d\frac{\pa f^\alpha}{\pa z^i}*dz^i\\
&=-*\frac{\pa^2 f^\alpha}{\pa z^i\pa z^j}dz^j\wedge * dz^i
-\frac{\pa f^\alpha}{\pa z^i}*d*dz^i
\end{align*}

We have
\[
d z^k\wedge *dz^j=0
\]
for all $k,j$ by the type considerations.

\begin{lemma}
\[
\sum_j\frac{\pa f^\alpha}{\pa z^j}*d*d z^j=
-\sum_j\frac{\pa f^\alpha}{\pa z^j}\Delta z^j=0
\]
\end{lemma}

{\bf Proof:}
Suppose that $z^j=y^j+\sqrt{-1} y^{j+N_1}$. Here $(y^1,\cdots, y^{2N_1})$
is
the  real local
coordinate. Suppose 
\[
g=g_{i\bar j}dz^i\otimes d\bar{z}^j
\]
is the invariant  Hermitian metric on $D$. The corresponding
Riemannian
metric is then
\begin{align*}
& G=G_{ab}dy^{a}dy^{b}\\
&=2(Re\,g_{i,\bar{j}}dy^idy^j
+Re\,g_{i,\bar{j}}dy^{i+N_1}dy^{j+N_1}\\
&-Im \,g_{i,\bar{j}}dy^{i+N_1}dy^j
+Im \,g_{i,\bar{j}}dy^{i}dy^{j+N_1})
\end{align*}

Here we assume the sum over $i,j,k,\cdots$ are from $1$ to $N_1$ and the sum over
$a,b,\cdots$ are from $1$ to $2N_1$. Thus
\begin{align*}
&\Delta z^k=\frac{1}{\sqrt{\det\, G}}\frac{\pa}{\pa x^{a}}
(\sqrt{\det\, G}G^{ab}\frac{\pa}{\pa x^b}z^k)\\
&=\frac {1}{\det g}\frac{\pa}{\pa x^a}
(\det g(G^{ak}+\sqrt{-1}G^{a(k+N_1)}))
\end{align*}

Now 
\begin{align*}
&(G^{ak}+\sqrt{-1}G^{a(k+N_1)})g_{k\bar j}\\
&=\frac 12 (G^{ak}+\sqrt{-1}G^{a(k+n)})(G_{kj}+\sqrt{-1}
G_{k,j+n})=\delta_{aj}+\sqrt{-1}\delta_{a(j+N_1)}
\end{align*}

Thus

\[
G^{ak}+\sqrt{-1}G^{a(k+n)}=
\left\{
\begin{array}{ll}
g^{k\bar a} & a\leq N_1\\
\sqrt{-1}g^{s,\bar{a-N_1}} & a>N_1
\end{array}
\right.
\]

So finally we have
\begin{align*}
&\Delta z^k=\frac{1}{\det g}\frac{\pa}{\pa x^l}
(\det g g^{k\bar l})
+\frac{1}{\det g}\frac{\pa}{\pa x^{l+N_1}}(\det g\sqrt{-1}g^{k\bar l})\\
&=2\frac{1}{\det g}\frac{\pa}{\pa\bar{z}^l}
(\det g g^{k\bar l})
=2g^{k\bar l} g^{i\bar j}(d{\omega})_{i\bar j\bar l}
\end{align*}
where $\omega=g_{i\bar j}d z^i\wedge d\bar{z}^j$ is the K\"ahler form of
the invariant Hermitian metric $g$ on $D$.

So we have
\[
\sum_j\frac{\pa f^\alpha}{\pa z^j}*d*d z^j=C\frac{\pa f^\alpha}{\pa z^j}
g^{j\bar{k}}g^{r\bar{s}}
(d\omega)_{r\bar{s}\bar{k}}
\]
where $C$ is a constant. Thus we have to prove
\[
C\frac{\pa f^\alpha}{\pa z^j} g^{j\bar{k}}g^{r\bar{s}}
(d\omega)_{r\bar{s}\bar{k}}=0
\]

{\bf Claim:}
\[
C\frac{\pa f^\alpha}{\pa z^j} g^{j\bar{k}}g^{r\bar{s}}
(d\omega)_{r\bar{s}\bar{k}}
=C<df^\alpha\wedge\omega,\bar{d\omega}>
\]

{\bf Proof:} By the definition of the inner product of the differential
forms.

\qed

Now let $x\in D$. Then we have the symmetry $\sigma_x$ defined in 
Equation~\ref{sym}. Now since $f$ is constant along each fiber, there is
an 
$\tilde{f}: G/K\rightarrow G/K$ such that 
\[
f=\tilde{f}\circ p
\]
Using $\sigma_x p=p\sigma_x$ and
$\sigma_x d\tilde f=-d\tilde f$, we have
\[
\sigma_x df^\alpha=\sigma_x d\tilde f^\alpha p=p\sigma_x d\tilde f^\alpha
=-p d\tilde f^\alpha =-df^\alpha
\]
On the other hand 
$\sigma_x\omega=\omega$ in Lemma~\ref{lm41}.
Thus
\[
<df^\alpha\wedge\omega, \bar{d\omega}>=\sigma_x<df^\alpha\wedge\omega,
\bar{d\omega}>
=-<df^\alpha\wedge\omega, \bar{d\omega}>
\]
The lemma is proved.\qed

{\bf Continuation of the Proof:} Now we know $f$ is a harmonic map with respect to the
metric $g$. 
In particular. $p$ is
harmonic map with respect to $g$. Since $p$ is a surjective map, by the
uniqueness theorem of the Harmonic map, we know $p=f$. But this
is a contradiction, because $f$ is a holomorphic map but $p$ is not. So
we have proved that
there are no K\"ahler metrics on $D$ which are $\Gamma$-invariant.

\bibliographystyle{acm}
\bibliography{bib}

\begin{thebibliography}{10}

\bibitem{Bl}
{\sc Black, M.}
\newblock {Harmonic Maps into Homogeneous Spaces}.
\newblock In {\em Pitman Research Notes in Mathematics Series\/} (1991),
  Longman Scientific and Technical.

\bibitem{Br}
{\sc Bryant, R.}
\newblock {Lie Groups and Twister Spaces}.
\newblock {\em Duke. Math J. 52\/} (1985), 233--261.

\bibitem{BS}
{\sc Burstall, F., and Salamons, S.}
\newblock {Tournaments, Flags and Harmonic Maps}.
\newblock {\em Math. Ann 277\/} (1987), 249--265.

\bibitem{ES}
{\sc Eells, J., and Sampson, J.}
\newblock {Harmonic Mappings of Riemannian Manifolds}.
\newblock {\em Amer. Jour. Math 86\/} (1964), 109--160.

\bibitem{Fd}
{\sc Friedman, R.}
\newblock {On threefolds with trivial canonical bundle}.
\newblock In {\em Complex geometry and Lie theory (Sundance, UT, 1989)\/}
  (1992), S.-T. Yau, Ed., International Press, pp.~103--134.

\bibitem{Gr}
{\sc Griffiths, P.}, Ed.
\newblock {\em {Topics in Transcendental Algebraic Geometry}}, vol.~106 of {\em
  Ann. Math Studies}.
\newblock Princeton University Press, 1984.

\bibitem{GS}
{\sc Griffiths, P., and Schmid, W.}
\newblock {Locally homogeneous complex manifolds}.
\newblock {\em Acta Math 123\/} (1969), 253--302.

\bibitem{Hel}
{\sc Helgason, S.}
\newblock {\em Differential Geometry, Lie Groups and Symmetric Spaces}.
\newblock Academic Press, 1978.

\bibitem{JY2}
{\sc Jost, J., and Yau, S.-T.}
\newblock {Harmonic maps and group representations}.
\newblock In {\em Differential geometry, Pitman Monographs Surveys Pure Appl.
  Math\/} (1991), Longman Sci. Tech., Harlow, pp.~241--259.

\bibitem{JY1}
{\sc Jost, J., and Yau, S.-T.}
\newblock {Harmonic mappings and algebraic varieties over function fields}.
\newblock {\em Amer. J. Math 115(6)\/} (1993), 1197--1227.

\bibitem{L}
{\sc Labourie, F.}
\newblock {Existence d'applications harmoniques tordues a valeurs dans les
  varietes a courbure negative}.
\newblock {\em Proc. Amer. Math. Soc 111(3)\/} (1991), 877--882.

\bibitem{Lu}
{\sc Lu, Z.}
\newblock {\em {On the Geometry of the Moduli Spaces}}.
\newblock PhD thesis, Courant Institute, New York University, New York, NY
  10012, June 1997.

\bibitem{M}
{\sc Murakami, S.}
\newblock {\em {Introduction to Homogeneous Spaces(in Chinese)}}.
\newblock Science and Technology Press in Shanghai, 1983.

\bibitem{Raj}
{\sc Rajan, C.~S.}
\newblock {Vanishing Theorems for Cohomologies of Automorphic Vector Bundles}.
\newblock {\em J. Diff. Geom. 43(2)\/} (1996), 376--409.

\bibitem{Siu1}
{\sc Siu, Y.~T.}
\newblock {The Complex Analyticity of Harmonic Maps and Strong Rigidity Compact
  K\"ahler Manifolds}.
\newblock {\em Ann. of Math. 112\/} (1980), 73--111.

\bibitem{T1}
{\sc Tian, G.}
\newblock {Smoothness of the Universal Deformation Space of Compact Calabi-Yau
  Manifolds and its Peterson-Weil Metric}.
\newblock In {\em Mathematical aspects of string theory\/} (1987), S.-T. Yau,
  Ed., vol.~1, World Scientific, pp.~629--646.

\bibitem{V2}
{\sc Viehweg, E.}
\newblock {\em {Quasi-Projective Moduli for Polarized Manifolds}}.
\newblock Ergebnisse der Mathematik und ihrer Grenzgebiete. Springer-Verlag,
  1991.

\bibitem{Yeung1}
{\sc Yeung, S.~K.}
\newblock {Compactification of K\"ahler Manifolds with Negative Ricci
  Curvature}.
\newblock {\em Invent. Math 106\/} (1991), 13--25.

\end{thebibliography}

\end{document}